\documentclass[fontsize=11pt,a4paper,reqno,numbers=noendperiod,bibliography=totoc]{scrartcl}
\usepackage[USenglish]{babel}
\usepackage[T1]{fontenc}
\usepackage[latin10]{inputenc}
\usepackage{fouriernc}
\usepackage{amsmath,amssymb,amsthm}
\usepackage{dsfont,upgreek}
\usepackage[obeyspaces]{url}
\usepackage{graphicx,tikz,pgfplots}
\usepackage{enumitem}
\usepackage{algorithm,algpseudocode}         
\usepackage{algorithm_lines}
\usepackage{cite} 

\AtBeginDocument{
\hypersetup{bookmarks=true,colorlinks=true,menucolor=yellow,citecolor=green!80!black,linkcolor=red!70!black,filecolor=magenta,urlcolor=magenta,breaklinks,pdfauthor={Bernardo Gonz\'alez Merino, Thomas Jahn, Alexandr Polyanskii and Gerd Wachsmuth},pdftitle={Hunting for reduced polytopes}}
}

\usetikzlibrary{calc}

\newcommand{\RR}{\mathds{R}}
\newcommand{\abs}[1]{\left\lvert#1\right\rvert}
\newcommand{\norm}[1]{\lVert#1\rVert}
\newcommand{\setn}[1]{\left\{#1\right\}}
\newcommand{\bigset}[1]{\bigl\{#1\bigr\}}
\newcommand{\bigp}[1]{\bigl(#1\bigr)}
\newcommand{\setcond}[2]{\left\{#1 \::\: #2\right\}}
\newcommand{\defeq}{\mathrel{\mathop:}=}
\newcommand{\skpr}[2]{\left\langle#1 \,\middle\vert\, #2\right\rangle}
\newcommand{\mc}[2]{B\!\left(#1,#2\right)}
\newcommand{\ms}[2]{S\!\left(#1,#2\right)}
\newcommand{\lr}[1]{\left(#1\right)}

\theoremstyle{plain}  
\newtheorem{Satz}{Theorem}[section]
\newtheorem{Lem}[Satz]{Lemma}
\newtheorem{Kor}[Satz]{Corollary}

\theoremstyle{definition} 
\newtheorem{Def}[Satz]{Definition}

\theoremstyle{remark}
\newtheorem*{myproof}{Proof of \hyperref[thm:nonreduced]{Theorem~\ref*{thm:nonreduced}}}

\DeclareMathOperator{\aff}{aff}
\DeclareMathOperator{\co}{co}
\DeclareMathOperator{\cone}{pos}
\DeclareMathOperator{\sbd}{sbd}

\let \subset \subseteq
\let \piup \uppi

\addtokomafont{caption}{\small}
\setkomafont{captionlabel}{\bfseries}

\newcommand{\mscLink}[1]{\href{http://www.ams.org/mathscinet/msc/msc2010.html?t=#1}{#1}}

\begin{document}
\parindent 0pt
\title{Hunting for reduced polytopes}
\author{Bernardo Gonz\'{a}lez Merino\footnote{TU Munich, Zentrum Mathematik, bg.merino@tum.de}
\thanks{The first author is partially supported by Consejer\'{i}a de Industria, Turismo, Empresa e Innovaci\'{o}n de la CARM through Fundaci\'{o}n S\'{e}neca, Agencia de Ciencia y Tecnolog\'{i}a de la Regi\'{o}n de Murcia, Programa de Formaci\'{o}n Postdoctoral de Personal Investigador 19769/PD/15 and project 19901/GERM/15, Programme in Support of Excellence Groups of the Regi\'{o}n de Murcia, and by MINECO project reference MTM2015-63699-P, Spain.}
	\and Thomas Jahn\footnote{TU Chemnitz, Faculty of Mathematics, thomas.jahn@mathematik.tu-chemnitz.de}
	\and Alexandr Polyanskii\footnote{Moscow Institute of Physics and Technology, Technion, Institute for Information Transmission Problems RAS, alexander.polyanskii@yandex.ru}
\thanks{The author was partially supported by the Russian Foundation for Basic Research, grants No. 15-31-20403 (mol\_a\_ved), No. 15-01-99563 A, No. 15-01-03530 A.}
	\and Gerd Wachsmuth\footnote{TU Chemnitz, Faculty of Mathematics, gerd.wachsmuth@mathematik.tu-chemnitz.de}
}
\date{\today}
\maketitle

\begin{abstract}
We show that there exist reduced polytopes in three-dimensional Euclidean space. This partially answers the question posed by Lassak \cite{Lassak1990} on the existence of reduced polytopes in $d$-dimensional Euclidean space for $d\geq 3$. Moreover, we prove a novel necessary condition on reduced polytopes in three-dimensional Euclidean space.
\end{abstract}


\textbf{MSC(2010):} \mscLink{52B10}

\section{Introduction}\label{chap:introduction}
Constant width bodies, i.e., convex bodies for which parallel supporting hyperplanes have constant distance, have a long and rich history in mathematics \cite{ChakerianGr1983}. Due to Meissner \cite{Meissner1911}, constant width bodies in Euclidean space can be characterized by \emph{diametrical completeness}, that is, the property of not being properly contained in a set of the same diameter. Constant width bodies also belong to a related class of \emph{reduced} convex bodies introduced by Heil \cite{Heil1978}. This means that constant width bodies do not properly contain a convex body of same minimum width. Remarkably, the classes of reduced bodies and constant width bodies do not coincide, as a regular triangle in the Euclidean plane shows.

Reduced bodies are extremal in remarkable inequalities for prescribed minimum width, as in Steinhagen's inequality \cite{ChakerianGr1983} (minimum inradius), or others that surprisingly still remain unsolved, namely, P\'al's problem \cite{Pal1921} (minimum volume). While the regular simplex (and any of its reduced subsets) is extremal for Steinhagen's, it is extremal only in the planar case for P\'al's problem. The reason is that while the regular triangle is reduced, this is no longer the case for the regular simplex in $\RR^d$, $d\geq 3$. Indeed, Heil conjectured \cite{Heil1978} that a certain reduced subset of the regular simplex is extremal for P\'al's problem. 
Heil also observed that some reduced body has to be extreme for P\'al's problem when replacing
volume by quermassintegral. The existence of reduced polytopes, and the fact that
smooth reduced sets are of constant width (cf.~\cite{Heil1978}), opens the door to conjecture some of
them as minimizers. In full generality, any non-decreasing-inclusion functional of convex
bodies with prescribed minimum width, attains its minimum at some reduced body.
P\'{a}l's problem restricted to constant width sets is the well-known Blaschke--Lebesgue problem, cf.\ \cite{ChakerianGr1983}, solved only in the planar case, where the Reuleaux triangle is the minimizer of the area, and Meissner's bodies are conjectured to be extremal in the three-dimensional space, see \cite[pp.\ 216, 217]{LassakMa2011} for an extended discussion. 
Note that P\'{a}l's problem has also been investigated in other geometrical settings such as Minkowskian planes \cite{Averkov2005} or spherical geometry, cf.~\cite[pp.\ 96, 97]{Bezdek2013} and \cite{LassakMusielak2016}.

Reduced bodies in the Euclidean space have been extensively studied in \cite{Lassak1990,Lassak2003b,LassakMa2011}, and the concept of reducedness has been translated to finite-dimensional normed spaces \cite{LassakMa2005,LassakMa2014,Lassak2012}. In reference to the existence of reduced polygons in the Euclidean plane, Lassak \cite{Lassak1990} posed the question whether there exist reduced polytopes in Euclidean $d$-space for $d\geq 3$. Several authors addressed the search for reduced polytopes in finite-dimensional normed spaces \cite{Lassak2006,AverkovMa2008a,AverkovMa2008b,MartiniWe2002,MartiniSw2004b}. For Euclidean space starting from dimension $3$ several classes of polytopes such as
\begin{itemize}
\item{polytopes in $\RR^d$ with $d+2$ vertices, $d+2$ facets, or more vertices than facets \cite[Corollary~7]{AverkovMa2008b},}
\item{centrally symmetric polytopes \cite[Claim~2]{LassakMa2005},}
\item{simple polytopes, i.e., polytopes in $\RR^d$ where each vertex is incident to $d$ edges (like polytopal prisms, for instance) \cite[Corollary~8]{AverkovMa2008b},}
\item{pyramids with polytopal base \cite[Theorem~1]{AverkovMa2008a}, and in particular simplices \cite{MartiniWe2002,MartiniSw2004b}}
\end{itemize}
were proved to be \emph{not} reduced. The purpose of the present article is two-fold. After proving a novel necessary condition for reduced polytopes in three-dimensional Euclidean space in \hyperref[chap:nonreduced]{Section~\ref*{chap:nonreduced}}, we present a reduced polytope in $\RR^3$ in \hyperref[chap:reduced_polytope]{Section~\ref*{chap:reduced_polytope}}. The validity of our example can be checked using the algorithm provided in \hyperref[chap:algorithm]{Section~\ref*{chap:algorithm}}.

\section{Notation and basic results}
Throughout this paper, we work in $d$-dimensional Euclidean space, that is, the vector space $\RR^d$ equipped with the inner product $\skpr{x}{y}\defeq \sum_{i=1}^d x_iy_i$ and the norm $\norm{x}\defeq\sqrt{\skpr{x}{x}}$, where $x=(x_1,\ldots,x_d)$ and $y=(y_1,\ldots,y_d)$ denote two points in $\RR^d$. A subset $K\subset\RR^d$ is said to be \emph{convex} if the line segment
\begin{equation*}
[x,y]\defeq\setcond{\lambda x+(1-\lambda)y}{0\leq\lambda\leq 1}
\end{equation*}
is contained in $K$ for all choices of $x,y\in K$. Convex compact subsets of $\RR^d$ having non-empty interior are called \emph{convex bodies}. The smallest convex superset of $K\subset\RR^d$ is called its \emph{convex hull} $\co(K)$, whereas the smallest affine subspace of $\RR^d$ containing $K$ is denoted by $\aff(K)$, the \emph{affine hull} of $K$. The \emph{affine dimension} $\dim(K)$ of $K$ is the dimension of its affine hull. The \emph{support function} $h(K,\cdot):\RR^d\to\RR$ of $K$ is defined by
\begin{equation*}
h(K,u)\defeq\sup\setcond{\skpr{u}{x}}{x\in K}.
\end{equation*}
For $u\in\RR^d\setminus\setn{0}$, the hyperplane $H(K,u)\defeq\setcond{x\in\RR^d}{\skpr{u}{x}=h(K,u)}$ is a \emph{supporting hyperplane} of $K$. The \emph{width} of $K$ in direction $u\in\RR^d$, defined by
\begin{equation*}
w(K,u)\defeq h(K,-u)+h(K,u)
\end{equation*}
equals the distance of the supporting hyperplanes $H(K,\pm u)$ multiplied by $\norm{u}$. The \emph{minimum width} of $K$ is $\omega(K)\defeq\inf\setcond{w(K,u)}{\norm{u}=1}$. A \emph{polytope} is the convex hull of finitely many points. The boundary of a polytope consists of \emph{faces}, i.e., intersections of the polytope with its supporting hyperplanes. We shall refer to faces of affine dimension $0$, $1$, and $d-1$ as \emph{vertices}, \emph{edges}, and \emph{facets}, respectively.
Faces of polytopes are lower-dimensional polytopes and shall be denoted by the list of their vertices. (A face which is denoted in this way can be reconstructed by taking the convex hull of its vertices.)
By definition, attainment of the minimal width of a polytope $P$ is related to a binary relation on faces of $P$ called \emph{strict antipodality}, see \cite{AverkovMa2008a}.
\begin{Def}\label{def:antipodal}
Let $P\subset\RR^d$ be a polytope. Distinct faces $F_1$, $F_2$ of $P$ are said to be \emph{strictly antipodal} if there exists a direction $u\in\RR^d$, $\norm{u}=1$, such that $H(P,u)\cap P=F_1$ and $H(P,-u)\cap P=F_2$.
\end{Def}

Gritzmann and Klee \cite[(1.9)]{GritzmannKl1992} formulated a necessary condition on strictly antipodal pairs whose distance equals the minimum width. Here, $F_1+F_2=\setcond{x+y}{x\in F_1,y\in F_2}$ denotes the Minkowski sum of sets $F_1,F_2\subset\RR^d$. (The set $F_1-F_2$ is defined analogously, $F_1\pm v$ shall be used as an abbreviation for $F_1\pm\setn{v}$ whenever $v\in\RR^d$, and, using the above conventions, $v_1v_2-v_3v_4=[v_1,v_2]-[v_3,v_4]$ for $v_1,\ldots,v_4\in\RR^d$.)
\begin{Satz}
\label{thm:width_polytope}
Suppose that $P\subset\RR^d$ is a polytope with non-empty interior, and that $F_1$ and $F_2$ are a strictly antipodal pair of faces of $P$ whose distance is equal to $\omega(P)$.
Then,
\begin{equation}
	\label{eq:gritzmann_klee}
\dim(F_1+F_2)=d-1,
\end{equation}
with $\dim(F_1)=\dim(F_2)=d-1$ when $P$ is centrally symmetric.
\end{Satz}
For arbitrary subsets $A,B\subset\RR^d$, we shall denote by 
\begin{equation*}
\rho(A,B)=\inf\setcond{\norm{x-y}}{x\in\aff(A),y\in\aff(B)}
\end{equation*}
the minimal distance between points of $\aff(A)$ and $\aff(B)$. In the situation of \hyperref[thm:width_polytope]{Theorem~\ref*{thm:width_polytope}}, $\rho(F_1,F_2)$ is then said distance between the respective parallel supporting hyperplanes of $P$.

The following definition by Heil \cite{Heil1978} is central to the present investigation.
\begin{Def}
A convex body $K$ is said to be \emph{reduced} if we have $\omega(K^\prime)<\omega(K)$ for all convex bodies $K^\prime\subsetneq K$.
\end{Def}
Reduced polytopes can be characterized using vertex-facet distances, see \cite[Theorem~4]{AverkovMa2008b} and \cite[Theorem~1]{Lassak2006} for the following result.

\begin{Satz}
\label{thm:polytope_reduced}
A polytope $P\subset\RR^d$ is reduced if and only if for every vertex $v$ of $P$, there exists a strictly antipodal facet $F$ of $P$ such that the distance between $v$ and $\aff(F)$ equals $\omega(P)$.
\end{Satz}

Strongly related, there is also the following necessary condition on the orthogonal projection of a vertex onto one of its strictly antipodal facets at the correct distance, see \cite[Lemma~2]{AverkovMa2008a}.
\begin{Satz}\label{thm:reduced_implication}
Assume that $P\subset\RR^d$ is a reduced polytope. Then for every vertex $v$ of $P$ there exists a facet $F$ of $P$ such that $\setn{v}$ is strictly antipodal to $F$, the orthogonal projection $w$ of $v$ onto $\aff(F)$ lies in the relative interior of $F$, and the distance from $v$ to $w$ is equal to $\omega(P)$.
\end{Satz}

\section{A class of non-reduced polytopes}\label{chap:nonreduced}
In this section, we prove the following necessary condition on reduced polytopes in three-dimensional Euclidean space. 

\begin{Satz}\label{thm:nonreduced}
Suppose that $P\subset\RR^3$ is a reduced polytope. Let $F$ be a facet of $P$ with edges $a_1a_2,\ldots,a_{k-1}a_k, a_ka_1$, and let $v$ be a vertex of $P$. Suppose that in this clockwise order, $vv_{1,1},\ldots,vv_{1,i_1}$, $vv_{2,1},\ldots,vv_{2,i_2},\ldots,$ $vv_{k,1},\ldots,vv_{k,i_k}$, where $k,i_1,\ldots,i_k$ denote positive integers, are the edges incident to $v$. For any $j\in\setn{1,\ldots,k}$ and $l\in\setn{1,\ldots,i_j-1}$, let $F_{j,l}$ be the facet incident to $vv_{j,l}$ and $vv_{j,l+1}$. For $j\in\setn{1,\ldots,k-1}$, let $F_{j,i_j}$ be the facet incident to $vv_{j,i_j}$ and $vv_{j+1,1}$. Finally, denote by $F_{k,i_k}$ be the facet incident to $vv_{k,i_k}$ and $vv_{1,1}$. Then the following conditions cannot be true at the same time:
\begin{enumerate}[label={(\alph*)},leftmargin=*,align=left,noitemsep]
\item{The facet $F$ and the vertex $v$ are strictly antipodal, and $\rho(v,F)=\omega(P)$.\label{condition1}}
\item{For any $j\in\setn{1,\ldots,k}$, the edges $vv_{j,1}$ and $a_{j-1}a_j$ are strictly antipodal. (Take $a_0=a_k$.)\label{condition2}}
\item{For any $j\in\setn{1,\ldots,k}$, the facets $F_{j,1},\ldots,F_{j,i_j}$ are strictly antipodal to $a_j$. Moreover, there is a number $l\in\setn{1,\ldots,i_j}$ such that $\rho(a_j,F_{j,l})=\omega(P)$.\label{condition3}}
\end{enumerate}
\end{Satz}
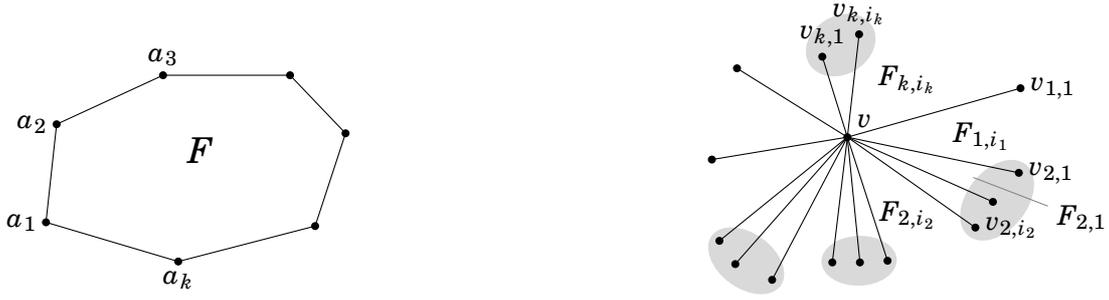
\begin{figure}[ht!]
\begin{center}
\begin{tikzpicture}[line cap=round,line join=round,>=stealth,x=1.0cm,y=1.0cm]
\coordinate (A1) at (0.46,1.05);
\coordinate (A2) at (0.6,2.35);
\coordinate (A3) at (2,3);
\coordinate (A4) at (3.67,3);
\coordinate (A5) at (4.4,2.23);
\coordinate (A6) at (4,1);
\coordinate (A7) at (2.2,0.53);

\draw (A1)--(A2)--(A3)--(A4)--(A5)--(A6)--(A7)--cycle;
\fill [color=black] (A1) circle (1.5pt) node[left]{$a_1$};
\fill [color=black] (A2) circle (1.5pt) node[left]{$a_2$};
\fill [color=black] (A3) circle (1.5pt) node[above]{$a_3$};
\fill [color=black] (A4) circle (1.5pt);
\fill [color=black] (A5) circle (1.5pt);
\fill [color=black] (A6) circle (1.5pt);
\fill [color=black] (A7) circle (1.5pt) node[below]{$a_k$};
\draw (2.5,2) node{\Large $F$};
\end{tikzpicture}
\hfill
\begin{tikzpicture}[line cap=round,line join=round,>=stealth,x=1.1cm,y=1.1cm]
\coordinate (V) at (2.28,5.81);
\coordinate (V11) at (4.35,6.4);
\coordinate (V21) at (4.33,5.38);
\coordinate (V22) at (4.02,5.03);
\coordinate (V23) at (3.81,4.72);
\coordinate (V31) at (2.76,4.32);
\coordinate (V32) at (2.43,4.3);
\coordinate (V33) at (2.1,4.3);
\coordinate (V41) at (1.38,4.09);
\coordinate (V42) at (0.94,4.28);
\coordinate (V43) at (0.75,4.56);
\coordinate (V51) at (0.66,5.54);
\coordinate (V61) at (0.96,6.64);
\coordinate (V71) at (1.98,6.78);
\coordinate (V72) at (2.42,7.05);

\fill [color=black!15,rotate around={51.84:(4.07,5.05)}] (4.07,5.05) ellipse (0.55 and 0.36);
\fill [color=black!15,rotate around={1.83:(2.42,4.32)}] (2.42,4.32) ellipse (0.45 and 0.3);
\fill [color=black!15,rotate around={-36.53:(1.07,4.32)}] (1.07,4.32) ellipse (0.51 and 0.33);
\fill [color=black!15,rotate around={31.83:(2.2,6.92)}] (2.2,6.92) ellipse (0.44 and 0.35);
\draw (V)--(V11);
\draw (V)--(V21);
\draw (V)--(V22);
\draw (V)--(V23);
\draw (V)--(V31);
\draw (V)--(V32);
\draw (V)--(V33);
\draw (V)--(V41);
\draw (V)--(V42);
\draw (V)--(V43);
\draw (V)--(V51);
\draw (V)--(V61);
\draw (V)--(V71);
\draw (V)--(V72);
\fill [color=black] (V) circle (1.5pt) node[above right]{$v$};
\fill [color=black] (V11) circle (1.5pt) node[right]{$v_{1,1}$};
\fill [color=black] (V21) circle (1.5pt) node[right]{$v_{2,1}$};
\fill [color=black] (V22) circle (1.5pt);
\fill [color=black] (3.81,4.72) circle (1.5pt) node[right]{$v_{2,i_2}$};
\fill [color=black] (2.76,4.32) circle (1.5pt);
\fill [color=black] (2.43,4.3) circle (1.5pt);
\fill [color=black] (2.1,4.3) circle (1.5pt);
\fill [color=black] (1.38,4.09) circle (1.5pt);
\fill [color=black] (0.94,4.28) circle (1.5pt);
\fill [color=black] (0.75,4.56) circle (1.5pt);
\fill [color=black] (0.66,5.54) circle (1.5pt);
\fill [color=black] (0.96,6.64) circle (1.5pt);
\fill [color=black] (1.98,6.78) circle (1.5pt) node[above]{$v_{k,1}$};
\fill [color=black] (2.42,7.05) circle (1.5pt) node[above]{$v_{k,i_k}$};
\draw (3.88,5.82) node{$F_{1,i_1}$};
\draw (3,4.87) node{$F_{2,i_2}$};
\fill (3.64,5.38) node[circle,pin={[pin distance=1cm]345:{$F_{2,1}$}}]{};
\draw (3,6.5) node{$F_{k,i_k}$};
\end{tikzpicture}
\end{center}\caption{Notation of \hyperref[thm:nonreduced]{Theorem~\ref*{thm:nonreduced}}: faces of $P$.}
\end{figure}
We prepare the proof of \hyperref[thm:nonreduced]{Theorem~\ref*{thm:nonreduced}} by three lemmas which rely on the geometry of the counterparts of convex polygons in spherical geometry. In order to avoid ambiguity, we fix the required notions and notation. The ball and the sphere with center $x\in\RR^3$ and radius $\alpha>0$ shall be denoted by
\begin{align*}
\mc{x}{\alpha}&\defeq\setcond{y\in\RR^3}{\norm{y-x}\leq\alpha},\\
\ms{x}{\alpha}&\defeq\setcond{y\in\RR^3}{\norm{y-x}=\alpha},
\end{align*}
respectively. The open half-space $H^<_{u,\alpha}\defeq\setcond{y\in\RR^3}{\skpr{u}{y}<\alpha}$ and its closed counterpart $H^\leq_{u,\alpha}\defeq\setcond{y\in\RR^3}{\skpr{u}{y}\leq\alpha}$ are bounded by the hyperplane $H^=_{u,\alpha}\defeq\setcond{y\in\RR^3}{\skpr{u}{y}=\alpha}$, where $u\in\ms{0}{1}$ denotes the outer normal unit vector of these three sets. Now fix a sphere $S=\ms{x_0}{\alpha}$. A \emph{great circle} of $S$ is the intersection of $S$ and a hyperplane $H^=_{u,\skpr{u}{x_0}}$. A \emph{hemisphere} of $S$ is the intersection of $S$ and an open half-space $H^<_{u,\skpr{u}{x_0}}$. Let $x,y\in S$ be contained in a hemisphere of $S$. There is exactly one great circle, denoted by $D_{x,y}$, passing through $x$ and $y$, and it is divided into two connected components by $x$ and $y$, one of which lies in the same hemisphere like $x$ and $y$. This connected component is called the \emph{arc} $\widearc{xy}$ whose length shall be denoted by $\abs{xy}$.

A \emph{cap} $C$ of $S$ is the intersection of $S$ and a closed half-space $H^\leq_{u,\beta}$ with $-\alpha<\beta<0$. The \emph{(spherical) boundary} $\sbd(C)$ of a cap $C=S\cap H^\leq_{u,\beta}$ is the circle $S\cap H^=_{u,\beta}$. The \emph{center} of the cap $C=S\cap H^\leq_{u,\beta}$ is the singleton $S\cap H^=_{u,-\alpha}$. A subset $A$ of $S$ which is contained in a hemisphere is said to be \emph{spherically convex} if for every choice of $x,y\in A$, the arc $\widearc{xy}$ is fully contained in $A$. Equivalently, a subset $A$ of $S$ is spherically convex if and only if the \emph{positive hull} $\cone(A-x_0)\defeq\setcond{\lambda u}{\lambda\geq 0, u\in A-x_0}$ is a convex set. A \emph{spherical polygon} is then the smallest spherically convex set containing a given finite subset of $S$. An arc $\widearc{xy}\subset S$ \emph{touches} a cap $C$ if $\aff\setn{x_0,x,y}\cap C$ is a singleton and is contained in $\widearc{xy}$.

\begin{Lem}\label{lem:cap_tangents}
Let $S^\prime$ be a hemisphere of $S$, $C\subset S^\prime$ be a cap, $x\in S^\prime\setminus C$. Furthermore, let $y_1,y_2\in\sbd(C)$ such that the arcs $\widearc{xy_1}$ and $\widearc{xy_2}$ touch $C$. Then $\abs{xy_1}=\abs{xy_2}$.
\end{Lem}
\begin{figure}[ht!]
\begin{center}
\begin{tikzpicture}[line cap=round,line join=round,>=stealth,x=0.5cm,y=0.5cm]
\draw [color=black] plot[domain=90:135,smooth,variable=\t]({-5*cos(60)*sin(-36.87)*sin(\t)-5*sin(60)*cos(\t)},{-5*sin(60)*sin(-36.87)*sin(\t)+5*cos(60)*cos(\t)});
\draw [color=black] plot[domain=90:135,smooth,variable=\t]({-5*cos(60)*sin(-36.87)*sin(\t)-5*sin(60)*cos(\t)},{5*sin(60)*sin(-36.87)*sin(\t)-5*cos(60)*cos(\t)});
\draw [color=black,line width=1pt] (0,0) circle(3);
\fill (4.16,0) circle(1.5pt) node[right]{$x$};
\fill (1.5,2.6) circle(1.5pt) node[below left]{$y_1$};
\fill (1.5,-2.6) circle(1.5pt) node[above left]{$y_2$};
\draw (-1.5,-1.9) node{\Large $C$};
\end{tikzpicture}
\end{center}\caption{Illustration of \hyperref[lem:cap_tangents]{Lemma~\ref*{lem:cap_tangents}}.}
\end{figure}
The proof of \hyperref[lem:cap_tangents]{Lemma~\ref*{lem:cap_tangents}} is left to the interested reader.
\begin{Lem}\label{lem:cap_tangents_2}
Let $S^\prime$ be a hemisphere of $S=\ms{x_0}{\alpha}$, $C\subset S^\prime$ be a cap, $a$, $b\in S^\prime\setminus C$. Assume that $\widearc{ab}\cap C$ contains at least two points. Furthermore, let $x,y\in\sbd(C)$ be such that $\widearc{ax}$ and $\widearc{by}$ touch $C$. Then $\abs{ab}\geq\abs{ax}+\abs{by}$.
\end{Lem}
\begin{figure}[ht!]
\begin{center}
\begin{tikzpicture}[line cap=round,line join=round,x=0.75cm,y=0.75cm]
\draw [rotate around={0:(0,-1.5)},line width=1pt] (0,-1.5) ellipse (2.3 and 2.16);
\draw [rotate around={-35:(0,0)}] plot[domain=61:98,smooth,variable=\t]({5*cos(\t)},{1*sin(\t)});
\draw [rotate around={35:(0,0)}] plot[domain=82:127,smooth,variable=\t]({5*cos(\t)},{1*sin(\t)});
\draw [rotate around={5:(0,0)}] plot[domain=234:298,smooth,variable=\t]({5*cos(\t)},{1*sin(\t)});

\draw(0,0) circle (1);
\fill [color=black] (-1.41,0.21) circle (1.5pt) node[left]{$x$};
\fill [color=black] (0.57,0.82) circle (1.5pt) node[above]{$a^\prime$};
\fill [color=black] (-0.57,0.82) circle (1.5pt) node[above]{$b^\prime$};
\fill [color=black] (0,-1.69) circle (1.5pt) node[right]{$z$};
\fill [color=black] (0,-1) circle (1.5pt) node[above]{$t^\prime$};
\fill [color=black] (0,1.21) circle (1.5pt) node[above]{$t$};
\fill [color=black] (1.41,0.21) circle (1.5pt) node[right]{$y$};
\fill [color=black] (-2.91,-1.06) circle (1.5pt) node[left] {$a$};
\fill [color=black] (2.48,-0.66) circle (1.5pt) node[right] {$b$};
\fill [color=black] (0,0) circle (1.5pt) node[right] {$z^\prime$};
\draw (-1.5,-2.55) node{\Large $C$};
\draw (-0.65,-0.25) node{$C^\prime$};
\end{tikzpicture}
\end{center}\caption{Illustration of \hyperref[lem:cap_tangents_2]{Lemma~\ref*{lem:cap_tangents_2}}.}
\end{figure}
\begin{proof}
Without loss of generality we assume that $x$ and $y$ lie in the same hemisphere bounded by $D_{a,b}$. Choose $t\in D_{a,x}\cap D_{b,y}$ such that that $x\in\widearc{at}$ and $y\in\widearc{bt}$, i.e., the points $x$, $y$, $t$ lie on the same hemisphere bounded by $D_{a,b}$. Consider the incircle $C^\prime$ of the spherical triangle with vertices $a$, $b$, and $t$, that is, the largest cap contained in this triangle. Denote its center by $z^\prime$. The incircle $C^\prime$ touches all sides of the spherical triangle; we thus set $\setn{a^\prime}\defeq C^\prime\cap\widearc{bt}$, $\setn{b^\prime}\defeq C^\prime\cap\widearc{at}$, and $\setn{t^\prime}\defeq C^\prime\cap\widearc{ab}$. Note that $b^\prime\in\widearc{xt}$. (Else, if $z$ denotes the center of the cap $C$, the arcs $\widearc{z^\prime b^\prime}$ and $\widearc{xz}$ intersect in a point $u$. Hence $\abs{ux}\leq\abs{zx}<\alpha\piup/2$. But $\widearc{ux}$ and $\widearc{ub^\prime}$ are both orthogonal to $D_{a,t}$, i.e., $\abs{ux}=\abs{ub^\prime}=\alpha\piup/2$, a contradiction.) Analogously, we have that $a^\prime\in\widearc{yt}$. Therefore,
\begin{equation*}
\abs{ax}+\abs{by}\leq\abs{ab^\prime}+\abs{ba^\prime}=\abs{at^\prime}+\abs{bt^\prime}=\abs{ab}.
\end{equation*}
Note that we have used \hyperref[lem:cap_tangents]{Lemma~\ref*{lem:cap_tangents}} in the first equality.
\end{proof}
\begin{Lem}\label{lem:spherical_geometry}
Let $p_1,\ldots,p_k$, $k\geq 4$, be (in this cyclic order) the vertices of a spherical polygon lying in a hemisphere $S^\prime$, and let $C\subset S^\prime$ be a cap overlapping with the interior of $P$. Assume that $\widearc{p_kp_1}$ touches $C$, and that $\widearc{p_1p_2}$ and $\widearc{p_{k-1}p_k}$ are disjoint to $C$ or touch $C$. Furthermore, assume that for some $i\in\setn{2,\ldots,k-2}$, $\widearc{p_ip_{i+1}}$ has non-empty intersection with $C$. Then
\begin{equation*}
\abs{p_kp_1}+\sum_{j=2}^{k-2}\abs{p_jp_{j+1}}\geq\abs{p_1p_2}+\abs{p_{k-1}p_k}
\end{equation*}
\end{Lem}
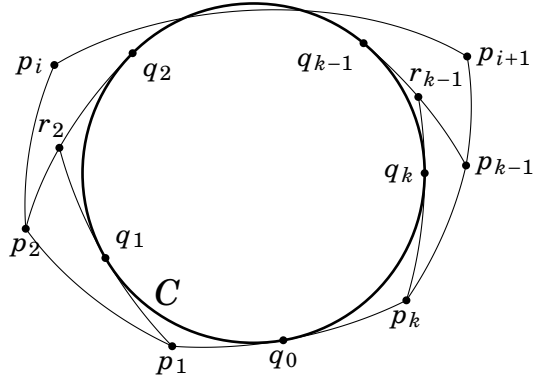
\begin{figure}[ht!]
\begin{center}
\begin{tikzpicture}[line cap=round,line join=round,x=0.75cm,y=0.75cm]
\draw[line width=1pt] (0,0) circle (3);

\draw [rotate around={90:(0,0)},color=black] plot[domain=243:285.5,smooth,variable=\t]({5*cos(\t)},{3*sin(\t)});
\draw [rotate around={10:(0,0)},color=black] plot[domain=247:297,smooth,variable=\t]({5*cos(\t)},{3*sin(\t)});
\draw [rotate around={140:(0,0)},color=black] plot[domain=237:270,smooth,variable=\t]({5*cos(\t)},{3*sin(\t)});
\draw [rotate around={-60:(0,0)},color=black] plot[domain=245:293,smooth,variable=\t]({5*cos(\t)},{3*sin(\t)});
\draw [rotate around={-135:(0,0)},color=black] plot[domain=270:315,smooth,variable=\t]({5*cos(\t)},{3*sin(\t)});

\draw [rotate around={61.213781292718814:(0,0)},color=black] plot[domain=262:293,smooth,variable=\t]({5*cos(\t)},{3.4633417115436624*sin(\t)});
\draw [rotate around={-25.39487568872483:(0,0)},color=black] plot[domain=230:270,smooth,variable=\t]({5*cos(\t)},{3.3764032295026056*sin(\t)});
\draw [rotate around={66.58626349307833:(0,0)},color=black] plot[domain=289:313,smooth,variable=\t]({5*cos(\t)},{3.552918316475135*sin(\t)});
\draw [rotate around={75.63200154251521:(0,0)},color=black] plot[domain=78.5:113,smooth,variable=\t]({5*cos(\t)},{3.936218153395842*sin(\t)});
\draw [rotate around={3.598070742230529:(0,0)},color=black] plot[domain=39.5:133,smooth,variable=\t]({5*cos(\t)},{2.8779544244285153*sin(\t)});

\fill (2.679563045280197,-2.2484203630108857) circle (1.5pt) node[below] {$p_k$};
\fill (3,0) circle (1.5pt) node[left] {$q_k$};
\fill (1.9283628290596173,2.2981333293569355) circle (1.5pt) node[below left] {$q_{k-1}$};
\fill (0.5209445330007902,-2.954423259036624) circle (1.5pt) node[below] {$q_0$};
\fill (-2.598076211353316,-1.5) circle (1.5pt) node[above right] {$q_1$};
\fill (-2.1213203435596424,2.1213203435596433) circle (1.5pt) node[below right] {$q_2$};
\fill (-3.405478863696411,0.4483397978867192) circle (1.5pt);
\fill (-3.55,0.8) node {$r_2$};
\fill (-1.4269480988832688,-3.060100073263126) circle (1.5pt) node[below] {$p_1$};
\fill (2.8890535823018912,1.3471878102475061) circle (1.5pt);
\fill (3.2,1.7) node {$r_{k-1}$};
\fill (3.725629019523878,0.135891646412324) circle (1.5pt) node[right] {$p_{k-1}$};
\fill (-3.995558560646749,-0.9823796197050645) circle (1.5pt) node[below] {$p_2$};
\fill (3.7453559935098006,2.063184164874058) circle (1.5pt) node[right] {$p_{i+1}$};
\fill (-3.4938011859457636,1.9121640163247595) circle (1.5pt) node[left] {$p_i$};
\draw (-1.5,-2.1) node{\Large $C$};
\end{tikzpicture}
\end{center}\caption{Illustration of \hyperref[lem:spherical_geometry]{Lemma~\ref*{lem:spherical_geometry}}.}
\end{figure}
\begin{proof}
The arc $\widearc{p_2p_{k-1}}$ has non-empty intersection with $C$. (Otherwise, $C$ is a subset of the spherical polygon with vertices $p_1$, $p_2$, $p_{k-1}$, and $p_k$, i.e., $C$ and $\widearc{p_ip_{i+1}}$ are disjoint.) Due to the triangle inequality,
\begin{equation*}
\sum_{j=2}^{k-2}\abs{p_jp_{j+1}}\geq\abs{p_2p_{k-1}}.
\end{equation*}
Thus, it remains to show that
\begin{equation*}
\abs{p_1p_k}+\abs{p_2p_{k-1}}\geq\abs{p_1p_2}+\abs{p_{k-1}p_k}.
\end{equation*}
Let $q_i\in\sbd(C)$ be such that the arc $\widearc{p_iq_i}$ touches $C$, $i\in\setn{1,2,k-1,k}$, $q_1,q_k\notin\widearc{p_1p_k}$. Let $q_2\in\sbd(C)$ be such that the arc $\widearc{p_2q_2}$ does not intersect $\widearc{p_1q_1}$. Analogously, let $q_{k-1}\in\sbd(C)$ be such that the arc $\widearc{p_{k-1}q_{k-1}}$ does not intersect $\widearc{p_kq_k}$. Let $q_0$ be the intersection point of $\sbd(C)$ and $\widearc{p_1p_k}$. The great circle through $p_1$ and $q_1$ intersects $\widearc{p_2q_2}$ in a point $r_2$. Similarly, the great circle through $p_k$ and $q_k$ intersects $\widearc{p_{k-1}q_{k-1}}$ in a point $r_{k-1}$.
Using the triangle inequality, \hyperref[lem:cap_tangents]{Lemma~\ref*{lem:cap_tangents}}, and \hyperref[lem:cap_tangents_2]{Lemma~\ref*{lem:cap_tangents_2}}, we obtain
\begin{align*}
\abs{p_1p_2}+\abs{p_kp_{k-1}}&\leq\abs{p_1q_1}+\abs{q_1r_2}+\abs{r_2p_2}+\abs{p_kq_k}+\abs{q_kr_{k-1}}+\abs{r_{k-1}p_{k-1}}\\
&=\abs{p_1q_0}+\abs{q_2r_2}+\abs{r_2p_2}+\abs{p_kq_0}+\abs{q_{k-1}r_{k-1}}+\abs{r_{k-1}p_{k-1}}\\
&=\abs{p_1p_k}+\abs{q_2p_2}+\abs{q_{k-1}p_{k-1}}\\
&\leq\abs{p_1p_k}+\abs{p_2p_{k-1}}.
\end{align*}
This completes the proof.
\end{proof}

\begin{myproof}
Suppose that $P$ is a reduced polytope which satisfies the conditions~\ref{condition1}, \ref{condition2}, and~\ref{condition3} of \hyperref[thm:nonreduced]{Theorem~\ref*{thm:nonreduced}}. Clearly, the facets of the difference polytope $P^\prime=P-P$ have the form $F_1-F_2$ with $F_1$ and $F_2$ being strictly antipodal faces of $P$ satisfying \hyperref[eq:gritzmann_klee]{Equation~(\ref*{eq:gritzmann_klee})}. Among the facets of the difference polytope $P^\prime=P-P$, there are $F_0=F-v$, $G_{j,l}=a_j-F_{j,l}$, $P_j=a_{j-1}a_j-vv_{j,1}$. Notice that $F_0$ and $G_{j,l}$ are congruent to $F$ and $F_{j,l}$, respectively, and that $P_j$ is a parallelogram. For $j\in\setn{1,\ldots,k}$ and $l\in\setn{1,\ldots,i_j}$, denote by $b_j=a_j-v$ the vertices of $F_0$, and let $c_{j,l}=a_j-v_{j,l}$, $d_j=a_j-v_{j+1,1}$. Obviously, $F_0$ and $P_j$ share a common edge $b_{j-1}b_j=a_{j-1}a_j-v$. 

\begin{figure}
\begin{center}
\begin{tikzpicture}[line cap=round,line join=round,>=stealth,x=1.0cm,y=1.0cm]
\coordinate (A1) at (0.46,1.05);
\coordinate (A2) at (0.6,2.35);
\coordinate (A3) at (2,3);
\coordinate (A4) at (3.67,3);
\coordinate (A5) at (4.4,2.23);
\coordinate (A6) at (4,1);
\coordinate (A7) at (2.2,0.53);
\coordinate (V) at (2.28,5.81);
\coordinate (V11) at (4.35,6.4);
\coordinate (V21) at (4.33,5.38);
\coordinate (V22) at (4.02,5.03);
\coordinate (V23) at (3.81,4.72);
\coordinate (V31) at (2.76,4.32);
\coordinate (V32) at (2.43,4.3);
\coordinate (V33) at (2.1,4.3);
\coordinate (V41) at (1.38,4.09);
\coordinate (V42) at (0.94,4.28);
\coordinate (V43) at (0.75,4.56);
\coordinate (V51) at (0.66,5.54);
\coordinate (V61) at (0.96,6.64);
\coordinate (V71) at (1.98,6.78);
\coordinate (V72) at (2.42,7.05);

\draw ($(A1)-(V)$)--($(A2)-(V)$)--($(A3)-(V)$)--($(A4)-(V)$)--($(A5)-(V)$)--($(A6)-(V)$)--($(A7)-(V)$)--cycle;
\fill [color=black] ($(A1)-(V)$) circle (1.5pt) node[above right]{$b_1$};
\fill [color=black] ($(A2)-(V)$) circle (1.5pt) node[below right]{$b_2$};
\fill [color=black] ($(A3)-(V)$) circle (1.5pt) node[below]{$b_3$};
\fill [color=black] ($(A4)-(V)$) circle (1.5pt);
\fill [color=black] ($(A5)-(V)$) circle (1.5pt);
\fill [color=black] ($(A6)-(V)$) circle (1.5pt);
\fill [color=black] ($(A7)-(V)$) circle (1.5pt) node[above]{$b_k$};
\fill [color=black] ($(A1)-(V11)$) circle (1.5pt) node[left]{$c_{1,1}$};
\fill [color=black] ($(A1)-(V21)$) circle (1.5pt) node[left]{$d_1$};
\fill [color=black] ($(A2)-(V21)$) circle (1.5pt) node[left]{$c_{2,1}$};
\fill [color=black] ($(A2)-(V22)$) circle (1.5pt);
\fill [color=black] ($(A2)-(V23)$) circle (1.5pt) node[above]{$c_{2,i_2}$};
\fill [color=black] ($(A2)-(V31)$) circle (1.5pt) node[above]{$d_2$};
\fill [color=black] ($(A3)-(V31)$) circle (1.5pt);
\fill [color=black] ($(A3)-(V32)$) circle (1.5pt);
\fill [color=black] ($(A3)-(V33)$) circle (1.5pt);
\fill [color=black] ($(A3)-(V41)$) circle (1.5pt);
\fill [color=black] ($(A4)-(V41)$) circle (1.5pt);
\fill [color=black] ($(A4)-(V42)$) circle (1.5pt);
\fill [color=black] ($(A4)-(V43)$) circle (1.5pt);
\fill [color=black] ($(A4)-(V51)$) circle (1.5pt);
\fill [color=black] ($(A5)-(V51)$) circle (1.5pt);
\fill [color=black] ($(A5)-(V61)$) circle (1.5pt);
\fill [color=black] ($(A6)-(V61)$) circle (1.5pt);
\fill [color=black] ($(A6)-(V71)$) circle (1.5pt);
\fill [color=black] ($(A7)-(V71)$) circle (1.5pt);
\fill [color=black] ($(A7)-(V72)$) circle (1.5pt);
\fill [color=black] ($(A7)-(V11)$) circle (1.5pt) node[below]{$d_k$};

\draw ($(A1)-(V)$)--($(A1)-(V11)$);
\draw ($(A1)-(V)$)--($(A1)-(V21)$);
\draw ($(A2)-(V)$)--($(A2)-(V21)$);
\draw ($(A2)-(V)$)--($(A2)-(V22)$);
\draw ($(A2)-(V)$)--($(A2)-(V23)$);
\draw ($(A2)-(V)$)--($(A2)-(V31)$);
\draw ($(A3)-(V)$)--($(A3)-(V31)$);
\draw ($(A3)-(V)$)--($(A3)-(V32)$);
\draw ($(A3)-(V)$)--($(A3)-(V33)$);
\draw ($(A3)-(V)$)--($(A3)-(V41)$);
\draw ($(A4)-(V)$)--($(A4)-(V41)$);
\draw ($(A4)-(V)$)--($(A4)-(V42)$);
\draw ($(A4)-(V)$)--($(A4)-(V43)$);
\draw ($(A4)-(V)$)--($(A4)-(V51)$);
\draw ($(A5)-(V)$)--($(A5)-(V51)$);
\draw ($(A5)-(V)$)--($(A5)-(V61)$);
\draw ($(A6)-(V)$)--($(A6)-(V61)$);
\draw ($(A6)-(V)$)--($(A6)-(V71)$);
\draw ($(A7)-(V)$)--($(A7)-(V71)$);
\draw ($(A7)-(V)$)--($(A7)-(V72)$);
\draw ($(A7)-(V)$)--($(A7)-(V11)$);
\draw ($(A1)-(V21)$)--($(A2)-(V21)$);
\draw ($(A2)-(V31)$)--($(A3)-(V31)$);
\draw ($(A3)-(V41)$)--($(A4)-(V41)$);
\draw ($(A4)-(V51)$)--($(A5)-(V51)$);
\draw ($(A5)-(V61)$)--($(A6)-(V61)$);
\draw ($(A6)-(V71)$)--($(A7)-(V71)$);
\draw ($(A7)-(V11)$)--($(A1)-(V11)$);

\draw ($(2.5,1.8)-(V)$) node{\Large $F_0$};
\draw ($(-0.7,1.9)-(V)$) node{\Large $P_2$};
\draw ($(0.6,0.5)-(V)$) node{\Large $P_1$};
\draw ($(3.3,0.2)-(V)$) node{\Large $P_k$};
\fill ($(-0.3,0.9)-(V)$) node[circle,pin={[pin distance=1cm]180:{$G_{1,i_1}$}}]{};
\fill ($(-0.3,2.7)-(V)$) node[circle,pin={[pin distance=1cm]170:{$G_{2,1}$}}]{};
\fill ($(0,3)-(V)$) node[circle,pin={[pin distance=1cm]100:{$G_{2,i_2}$}}]{};
\end{tikzpicture}
\end{center}\caption{Notation of \hyperref[thm:nonreduced]{Theorem~\ref*{thm:nonreduced}}: faces of $P^\prime$.}
\end{figure}

We have
\begin{align*}
\rho(0,F_0)&=\omega(P),\\
\rho(0,P_j)&\geq\omega(P)\qquad\text{for }j\in\setn{1,\ldots,k},\\
\text{and}\qquad\rho(0,G_{j,l})&\geq\omega(P)\qquad\text{for }j\in\setn{1,\ldots,k}\text{ and }l\in\setn{1,\ldots,i_j}.
\end{align*}
Moreover, for each $j\in\setn{1,\ldots,k}$, there is $l\in\setn{1,\ldots,i_j}$ such that $\rho(0,G_{j,l})=\omega(P)$.
Next, we prove the inequalities
\begin{equation}
\measuredangle b_{j-1}b_jb_{j+1}+\sum_{l=1}^{i_j-1}\measuredangle c_{j,l}b_jc_{j,l+1}+\measuredangle c_{j,i_j}b_jd_j\geq\measuredangle b_{j-1}b_jc_{j,1}+\measuredangle d_jb_jb_{j+1}\label{eq:angles}
\end{equation}
for $j\in\setn{1,\ldots,k}$, where we use the abbreviation $\measuredangle xyz=\arccos\lr{\frac{\skpr{x-y}{z-y}}{\norm{x-y}\norm{z-y}}}$ for the angle between vectors $x-y$ and $z-y$. Denote by $q$ the orthogonal projection of $0$ onto $\aff(F_0)$ and consider the ball $B=\mc{0}{\omega(P)}$ with center $0$ and radius $\omega(P)$ and the spheres $S_j=\ms{b_j}{\norm{b_j-q}}$, $j\in\setn{1,\ldots,k}$, with center $b_j$ and radius equal to the length of the segment $b_jq$. Let $C_j=S_j\cap B$ and $Q_j=S_j\cap(b_j+\cone(P^\prime-b_j))$. Then $Q_j$ is a spherical polygon contained in a hemisphere of $S_j$. 
(Choose a supporting hyperplane $H$ of $P^\prime$ with $H\cap P^\prime=\setn{b_j}$. Then $H\cap S_j$ is the spherical boundary of an admissible hemisphere.) Since
\begin{equation*}
C_j=S_j\cap B\subset B\subset P^\prime,
\end{equation*}
we also have 
\begin{equation*}
b_j+\cone(C_j-b_j)\subset b_j+\cone(P^\prime-b_j),
\end{equation*}
and thus
\begin{equation*}
C_j=S_j\cap (b_j+\cone(C_j-b_j))\subset S_j\cap (b_j+\cone(P^\prime-b_j))=Q_j.
\end{equation*}
Denote by $p_{j,1}$ the intersection point of $S_j$ and $b_j+\cone(b_{j-1}-b_j)$, by $p_{j,l}$, $l\in\setn{2,\dots, i_j+1}$, the intersection point of $S_j$ and $b_j+\cone(c_{j,l-1}-b_j)$, by $p_{j,i_j+2}$ the intersection point of $S_j$ and $b_j+\cone(d_j-b_j)$, and by $p_{j,i_j+3}$ the intersection point of $S_j$ and $b_j+\cone(b_{j+1}-b_j)$. Since $\rho(0,G_{j,l-1})=\omega(P)$ for some $l\in\setn{2,\ldots,i_j+1}$, the arc $\widearc{p_{j,l}p_{j,l+1}}$ touches the cap $C_j$. Analogously, $\widearc{p_{j,i_j+3}p_{j,1}}$ touches $C_j$. Notice that $\widearc{p_{j,1}p_{j,2}}$ and $\widearc{p_{j,i_j+2}p_{j,i_j+3}}$ intersect $C_j$ in at most one point each because $\rho(0,P_j)\geq\omega(P)$ and $\rho(0,P_{j+1})\geq\omega(P)$. Therefore, the assumptions from \hyperref[lem:spherical_geometry]{Lemma~\ref*{lem:spherical_geometry}} hold. In particular, since the lengths of arcs in \hyperref[lem:spherical_geometry]{Lemma~\ref*{lem:spherical_geometry}} correspond to angles in \hyperref[eq:angles]{Equation~(\ref*{eq:angles})}, the latter equation is true. Note that $\measuredangle d_jb_jb_{j+1}+\measuredangle b_jb_{j+1}c_{i+1,1}=\piup$ for all $j\in\setn{1,\ldots,k}$ because $P_j$ is a parallelogram. Adding the inequalities \hyperref[eq:angles]{(\ref*{eq:angles})} for $j\in\setn{1,\ldots,k}$ yields 
\begin{equation*}
\sum_{j=1}^k\lr{\measuredangle b_{j-1}b_jb_{j+1}+\sum_{l=1}^{i_j-1}\measuredangle c_{j,l}b_jc_{j,l+1}+\measuredangle c_{j,i_j}b_jd_j}\geq\sum_{j=1}^k\lr{\measuredangle b_{j-1}b_jc_{j,1}+\measuredangle d_jb_jb_{j+1}}=k\piup.
\end{equation*}
Hence,
\begin{equation}
\sum_{j=1}^k\lr{\measuredangle b_{j-1}b_jb_{j+1}+\sum_{l=1}^{i_j-1}\measuredangle c_{j,l}b_jc_{j,l+1}+\measuredangle c_{j,i_j}b_jd_j+\measuredangle b_{j-1}b_jc_{j,1}+\measuredangle d_jb_jb_{j+1}}\geq 2k\piup\label{eq:solid-angles}
\end{equation}
This last inequality contradicts the fact that the angles occuring on the left-hand side of \hyperref[eq:solid-angles]{Equation~(\ref*{eq:solid-angles})} are internal angles of the facets adjacent to $b_j$ at this vertex.
\end{myproof}

Note that the idea of the proof of \hyperref[thm:nonreduced]{Theorem~\ref*{thm:nonreduced}} is similar to Steinitz's approach in \cite{Steinitz1928}, where he constructed an example of non-circumscribable polytope in $\RR^3$. 

Specifying \hyperref[thm:nonreduced]{Theorem~\ref*{thm:nonreduced}} to $i_1=\ldots=i_k=1$, we obtain the following result.
\begin{Kor}\label{cor:nonreduced}
Suppose that $P\subset\RR^3$ is a reduced polytope. Let $F$ be a facet of $P$ with edges $a_1a_2,\ldots,a_{k-1}a_k, a_ka_1$, and let $v$ be a vertex of $P$. Suppose that in this clockwise order, $vv_1,\ldots,vv_2,\ldots,$ $vv_k$ are the edges incident to $v$. For any $j\in\setn{1,\ldots,k-1}$, let $F_j$ be the facet incident to $vv_j$ and $vv_{j+1}$. Finally, denote by $F_k$ be the facet incident to $vv_k$ and $vv_1$. Then the following conditions cannot be true at the same time:
\begin{enumerate}[label={(\alph*)},leftmargin=*,align=left,noitemsep]
\item{The facet $F$ and the vertex $v$ are strictly antipodal, and $\rho(v,F)=\omega(P)$.\label{first}}
\item{For any $j\in\setn{1,\ldots,k}$, the edges $vv_j$ and $a_{j-1}a_j$ are strictly antipodal. (Take $a_0=a_k$.)\label{second}}
\item{For any $j\in\setn{1,\ldots,k}$, the facet $F_j$ is strictly antipodal to $a_j$ and $\rho(a_j,F_j)=\omega(P)$.\label{third}}
\end{enumerate}
\end{Kor}
Now assume that $P\subset\RR^3$ is a combinatorially self-dual polytope, i.e., there exists an inclusion-reversing bijective map $\phi$ from the face lattice of $P$ onto itself, and for each face $F$ of $P$, $\phi(F)$ is the unique antipodal face of $F$. Then the conditions \ref{first}, \ref{second}, and \ref{third} in \hyperref[cor:nonreduced]{Corollary~\ref*{cor:nonreduced}} are satisfied at each vertex of $P$ which renders the polytope non-reduced. In particular, this applies to the case of pyramids in three-dimensional Euclidean space which belong to several classes of non-reduced polytopes mentioned in \hyperref[chap:introduction]{Section~\ref*{chap:introduction}}.

\section{A reduced polytope}\label{chap:reduced_polytope}
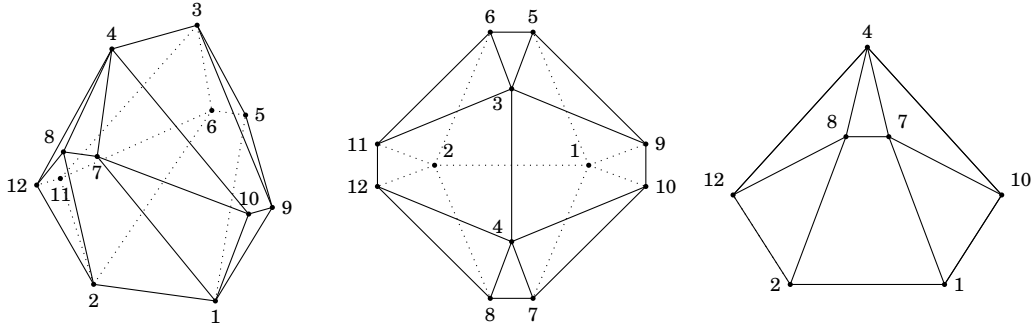
\begin{figure}
\begin{center}
\begin{tikzpicture}[
		dot/.style={fill,shape=circle,scale=0.2},
		every label/.style={scale=0.7}
]
\begin{axis}[%
width=1.8in,
height=1.8in,
view={35}{10},
scale only axis,
hide axis,
xmin=-0.6,
xmax= 0.6,
ymin=-0.6,
ymax= 0.6,
zmin=-0.6,
zmax= 0.6,
]

\pgfmathsetmacro{\t}{0.550000000000000}
\pgfmathsetmacro{\x}{0.617649095979951}
\pgfmathsetmacro{\s}{0.135138493102611}
\pgfmathsetmacro{\h}{0.098430025240916}
\pgfmathsetmacro{\r}{0.354718358670911}

\coordinate (A) at (axis cs: \r,   0, -\t);
\coordinate (B) at (axis cs:-\r,   0, -\t);
\coordinate (C) at (axis cs:  0,  \r,  \t);
\coordinate (D) at (axis cs:  0, -\r,  \t);
\coordinate (E) at (axis cs: \h,  \x,  \s);
\coordinate (F) at (axis cs:-\h,  \x,  \s);
\coordinate (G) at (axis cs: \h, -\x,  \s);
\coordinate (H) at (axis cs:-\h, -\x,  \s);
\coordinate (I) at (axis cs: \x,  \h, -\s);
\coordinate (J) at (axis cs: \x, -\h, -\s);
\coordinate (K) at (axis cs:-\x,  \h, -\s);
\coordinate (L) at (axis cs:-\x, -\h, -\s);

\node [label={below     : 1},dot] at (A) {};
\node [label={below     : 2},dot] at (B) {};
\node [label={above     : 3},dot] at (C) {};
\node [label={above     : 4},dot] at (D) {};
\node [label={right     : 5},dot] at (E) {};
\node [label={below     : 6},dot] at (F) {};
\node [label={below     : 7},dot] at (G) {};
\node [label={above left: 8},dot] at (H) {};
\node [label={right     : 9},dot] at (I) {};
\node [label={above     :10},dot] at (J) {};
\node [label={below     :11},dot] at (K) {};
\node [label={left      :12},dot] at (L) {};

\foreach \a/\b in {A/B,C/D,A/I,A/J,C/I,D/J,I/J,B/L,D/L,G/H,D/G,D/H,H/L,G/J,E/I,C/E,A/G,B/H} {
\edef\temp{\noexpand \draw (\a)--(\b);}
\temp
}

\foreach \a/\b in {B/K,E/F,C/F,F/K,A/E,B/F,C/K,K/L} {
\edef\temp{\noexpand \draw[dotted] (\a)--(\b);}
\temp
}
\end{axis}
\end{tikzpicture}
\begin{tikzpicture}[
		dot/.style={fill,shape=circle,scale=0.2},
		every label/.style={scale=0.7}
]
\begin{axis}[%
width=1.8in,
height=1.8in,
view={0}{90},
scale only axis,
hide axis,
xmin=-0.8,
xmax= 0.8,
ymin=-0.8,
ymax= 0.8,
]

\pgfmathsetmacro{\t}{0.550000000000000}
\pgfmathsetmacro{\x}{0.617649095979951}
\pgfmathsetmacro{\s}{0.135138493102611}
\pgfmathsetmacro{\h}{0.098430025240916}
\pgfmathsetmacro{\r}{0.354718358670911}

\coordinate (A) at (axis cs: \r,   0, -\t);
\coordinate (B) at (axis cs:-\r,   0, -\t);
\coordinate (C) at (axis cs:  0,  \r,  \t);
\coordinate (D) at (axis cs:  0, -\r,  \t);
\coordinate (E) at (axis cs: \h,  \x,  \s);
\coordinate (F) at (axis cs:-\h,  \x,  \s);
\coordinate (G) at (axis cs: \h, -\x,  \s);
\coordinate (H) at (axis cs:-\h, -\x,  \s);
\coordinate (I) at (axis cs: \x,  \h, -\s);
\coordinate (J) at (axis cs: \x, -\h, -\s);
\coordinate (K) at (axis cs:-\x,  \h, -\s);
\coordinate (L) at (axis cs:-\x, -\h, -\s);

\node [label={above  left: 1},dot] at (A) {};
\node [label={above right: 2},dot] at (B) {};
\node [label={below  left: 3},dot] at (C) {};
\node [label={above  left: 4},dot] at (D) {};
\node [label={above      : 5},dot] at (E) {};
\node [label={above      : 6},dot] at (F) {};
\node [label={below      : 7},dot] at (G) {};
\node [label={below      : 8},dot] at (H) {};
\node [label={right      : 9},dot] at (I) {};
\node [label={right      :10},dot] at (J) {};
\node [label={left       :11},dot] at (K) {};
\node [label={left       :12},dot] at (L) {};

\foreach \a/\b in {C/D,C/I,D/J,I/J,D/L,G/H,D/G,D/H,H/L,G/J,C/K,K/L,E/I,C/E,C/F,F/K,E/F} {
\edef\temp{\noexpand \draw (\a)--(\b);}
\temp
}

\foreach \a/\b in {A/B,A/I,A/J,A/G,A/E,B/K,B/H,B/F,B/L} {
\edef\temp{\noexpand \draw[dotted] (\a)--(\b);}
\temp
}
\end{axis}
\end{tikzpicture}
\begin{tikzpicture}[
		dot/.style={fill,shape=circle,scale=0.2},
		every label/.style={scale=0.7}
]

\begin{axis}[%
width=1.8in,
height=1.8in,
view={0}{0},
scale only axis,
hide axis,
xmin=-0.8,
xmax= 0.8,
ymin=-0.8,
ymax= 0.8,
zmin=-0.8,
zmax= 0.8,
]

\pgfmathsetmacro{\t}{0.550000000000000}
\pgfmathsetmacro{\x}{0.617649095979951}
\pgfmathsetmacro{\s}{0.135138493102611}
\pgfmathsetmacro{\h}{0.098430025240916}
\pgfmathsetmacro{\r}{0.354718358670911}

\coordinate (A) at (axis cs: \r,   0, -\t);
\coordinate (B) at (axis cs:-\r,   0, -\t);
\coordinate (C) at (axis cs:  0,  \r,  \t);
\coordinate (D) at (axis cs:  0, -\r,  \t);
\coordinate (E) at (axis cs: \h,  \x,  \s);
\coordinate (F) at (axis cs:-\h,  \x,  \s);
\coordinate (G) at (axis cs: \h, -\x,  \s);
\coordinate (H) at (axis cs:-\h, -\x,  \s);
\coordinate (I) at (axis cs: \x,  \h, -\s);
\coordinate (J) at (axis cs: \x, -\h, -\s);
\coordinate (K) at (axis cs:-\x,  \h, -\s);
\coordinate (L) at (axis cs:-\x, -\h, -\s);

\node [label={right      : 1},dot] at (A) {};
\node [label={left       : 2},dot] at (B) {};
\node [label={above      : 4},dot] at (D) {};
\node [label={above right: 7},dot] at (G) {};
\node [label={above left : 8},dot] at (H) {};
\node [label={above right:10},dot] at (J) {};
\node [label={above left :12},dot] at (L) {};

\foreach \a/\b in {A/B,C/D,A/I,A/J,C/I,D/J,I/J,B/L,D/L,G/H,D/G,D/H,H/L,G/J,C/K,K/L,A/G,B/H} {
\edef\temp{\noexpand \draw (\a)--(\b);}
\temp
}

\foreach \a/\b in {B/K} {
\edef\temp{\noexpand \draw[dotted] (\a)--(\b);}
\temp
}
\end{axis}
\end{tikzpicture}
\end{center}\caption{A reduced polytope $P$ with vertex numbers: oblique view (left), top view (middle), front view (right)}
\label{fig:our_polytope}
\end{figure}

In contrast to the various classes of polytopes which are shown to be non-reduced in the literature and in \hyperref[chap:nonreduced]{Section~\ref*{chap:nonreduced}}, we present a reduced polytope $P$ now. Consider the points
\begin{align*}
	v_{ 1} &\defeq ( r,  0, -t), &
	v_{ 2} &\defeq (-r,  0, -t), &
	v_{ 3} &\defeq ( 0,  r,  t), &
	v_{ 4} &\defeq ( 0, -r,  t), \\
	v_{ 5} &\defeq ( h,  x,  s), &
	v_{ 6} &\defeq (-h,  x,  s), &
	v_{ 7} &\defeq ( h, -x,  s), &
	v_{ 8} &\defeq (-h, -x,  s), \\
	v_{ 9} &\defeq ( x,  h, -s), &
	v_{10} &\defeq ( x, -h, -s), &
	v_{11} &\defeq (-x,  h, -s), &
	v_{12} &\defeq (-x, -h, -s).
\end{align*}
For properly chosen parameters $t,x,s,h,r>0$ the points $v_1,\ldots,v_{12}$ are the vertices of our polytope $P$.
The combinatorial structure of our polytope is shown in \hyperref[fig:our_polytope]{Figure~\ref*{fig:our_polytope}}.
The polytope $P$ possesses the same symmetry as the Johnson solid $J_{84}$ (however, not the same combinatorial structure).
Hence, it is sufficient to control few
facet-vertex and edge-edge distances.
In fact, we are going to solve the equations
\begin{align*}
	\rho(v_1,\; v_{3} \, v_{11} \, v_{12} \, v_{4}) &= 1, &
	\rho(v_1 \, v_{2},\; v_{ 3} \, v_{ 4}         ) &= \delta_1, &
	\rho(v_1 \, v_{9},\; v_{ 4} \, v_{ 8}         ) &= \delta_3, \\
	\rho(v_5,\; v_{2} \, v_{ 8} \, v_{12}         ) &= 1, &
	\rho(v_1 \, v_{5},\; v_{ 4} \, v_{ 8}         ) &= \delta_2,
\end{align*}
with respect to $t,x,s,h,r$.
Here, $\delta_1, \delta_2, \delta_3 \geq 1$ are suitably chosen.
By introducing the normal vectors
\begin{align*}
	n_1 &\defeq (v_{11} - v_3) \times (v_{12} - v_3), &
	n_4 &\defeq (v_{ 1} - v_5) \times (v_{ 4} - v_8), \\
	n_2 &\defeq (v_{ 8} - v_2) \times (v_{12} - v_2), &
	n_5 &\defeq (v_{ 1} - v_9) \times (v_{ 4} - v_8), \\
	n_3 &\defeq (v_{ 1} - v_2) \times (v_{ 3} - v_4) = (0,0,4\,r^2),
\end{align*}
where $u \times w$ denotes the usual cross product of the vectors $u,w \in \RR^3$,
these equations can be rewritten as
\begin{align*}
	\skpr{n_1}{v_1 - v_3}^2 -               \norm{n_1}^2 &= 0, &
	\skpr{n_3}{v_3 - v_1}^2 - \delta_1^2 \, \norm{n_3}^2 &= 0, &
	\skpr{n_5}{v_1 - v_4}^2 - \delta_3^2 \, \norm{n_5}^2 &= 0, \\
	\skpr{n_2}{v_5 - v_2}^2 -               \norm{n_2}^2 &= 0, &
	\skpr{n_4}{v_1 - v_4}^2 - \delta_2^2 \, \norm{n_4}^2 &= 0.
\end{align*}
Now, it is easy to see that the third equation (counting in columns) is equivalent to $2 \, t = \delta_1$.
Moreover, it is tedious to check that
we can factor out $h^2$ in the first equation and $(h + r - x)^2$ in the fifth.
Hence, we are going solve the four equations
\begin{equation}
	\label{eq:equations}
	\left.
		\begin{aligned}
			h^{-2}           \, \big(\skpr{n_1}{v_1 - v_3}^2 - \phantom{\delta_1^2} \, \norm{n_1}^2 \big) &= 0, &
			                         \skpr{n_2}{v_5 - v_2}^2 - \phantom{\delta_1^2} \, \norm{n_2}^2       &= 0, \\
			(h + r - x)^{-2} \, \big(\skpr{n_5}{v_1 - v_4}^2 -          \delta_3^2  \, \norm{n_5}^2 \big) &= 0, &
			                         \skpr{n_4}{v_1 - v_4}^2 -          \delta_2^2  \, \norm{n_4}^2       &= 0,
		\end{aligned}
		\qquad
	\right\}
\end{equation}
under $t = \delta_1/2$ with respect to the remaining variables $(x,s,h,r)$.
Note that each left-hand side of the four equations in \hyperref[eq:equations]{(\ref*{eq:equations})}
are multivariate polynomials of degree at most $6$ in the four unknowns $(x,s,h,r)$.

Numerically, we used $\delta_1 = 1.1$, $\delta_2 = 1.003$ and $\delta_3 = 1.004$
and solved equations \hyperref[eq:equations]{(\ref*{eq:equations})} by Newton's method starting with
$(x_0,s_0,h_0,r_0) = (0.62, 0.13, 0.09, 0.35)$.
This results in
\begin{align*}
	t &= 0.55, &
	x &\approx 0.6176490959800, &
	s &\approx 0.1351384931026, \\
	&&
	h &\approx 0.0984300252409, &
	r &\approx 0.3547183586709
\end{align*}
and the numerical residuum in the four equations is below $10^{-15}$.
Using Kantorovich's theorem, see \cite[Theorem~XVIII.1.6]{KantorovichAkilov1982},
it is possible to prove that
equations \hyperref[eq:equations]{(\ref*{eq:equations})} possess
an exact root in the neighborhood of our numerical solution.

Using these parameters, we can check that
the remaining distances are
\begin{align*}
	\rho( v_1 \, v_9 \, v_{10} ,\; v_{11} \, v_{12} ) &\approx 1.0433929735637, &
	\rho( v_5 \, v_9           ,\; v_{ 8} \, v_{12} ) &\approx 1.0126888049628.
\end{align*}
Thus, the width of our polytope using these parameters is really $1.0$,
see \hyperref[thm:width_polytope]{Theorem~\ref*{thm:width_polytope}}.
Consequently, our polytope is reduced by \hyperref[thm:polytope_reduced]{Theorem~\ref*{thm:polytope_reduced}}.

Since the Jacobian of the (left-hand sides of) equations \hyperref[eq:equations]{(\ref*{eq:equations})} with respect to $(x, s, h, r)$ is invertible at our point of interest,
it follows from the implicit function theorem that we also obtain a solution for small changes of the parameters
$\delta_1, \delta_2$ and $\delta_3$.
Hence, we obtain a whole family of reduced polytopes possessing three degrees of freedom.

\section{Evaluating your catches}\label{chap:algorithm}
It is quite a delicate and tedious procedure to check
the reducedness of a given polytope $P \subset \RR^3$.
Hence, we present an algorithm
based on \hyperref[thm:width_polytope]{Theorems~\ref*{thm:width_polytope}} and \hyperref[thm:polytope_reduced]{\ref*{thm:polytope_reduced}}.
It consists of two steps:
\begin{enumerate}
	\item
		Compute the width of $P$, compare \hyperref[thm:width_polytope]{Theorem~\ref*{thm:width_polytope}}.
	\item
		Check whether each vertex has a strictly antipodal facet, compare \hyperref[thm:polytope_reduced]{Theorem~\ref*{thm:polytope_reduced}}.
\end{enumerate}
An implementation in pseudocode is given in \hyperref[alg:red]{Algorithm~\ref*{alg:red}}.
In step~4 of the algorithm, we denoted by $e_1 \times e_2$
a vector normal to the skew edges $e_1$ and $e_2$.
Using \hyperref[thm:width_polytope]{Theorems~\ref*{thm:width_polytope}} and \hyperref[thm:polytope_reduced]{\ref*{thm:polytope_reduced}},
it is easy to check its correctness.
A Matlab implementation is provided at zenodo, see \cite{GonzalezMerinoJahnWachsmuth2016:2}.
\begin{algorithm}[ht]
	\begin{algorithmic}[1]
		\LState \textbf{input:} polytope $P \subset \RR^3$
		\LState set $w \leftarrow +\infty$
		\ForAll{skew pairs of edges $e_1$ and $e_2$ of $P$}
			\LState set $w \leftarrow \min\bigset{w, w\bigp{e_1 \times e_2 / \norm{e_1 \times e_2}, P}}$
		\EndFor
		\LState unmark all vertices of $P$
		\ForAll{facets $F$ of $P$}
			\LState compute the strictly antipodal face $\hat F$
			\LState set $\hat w \leftarrow \rho(F, \hat F)$
			\If{$\hat w < w$}
				\LState unmark all vertices of $P$
				\LState set $w \leftarrow \hat{w}$
			\EndIf
			\If{$\hat{F}$ consists of a single vertex $v$ \textbf{and} $\hat{w}= w$}
				\LState mark vertex $v$
			\EndIf
		\EndFor
		\LState\Return Are all vertices of $P$ marked?
	\end{algorithmic}
	\caption{Algorithm for checking reducedness of polytopes in $\RR^3$.}
	\label{alg:red}
\end{algorithm}

\section{Concluding remarks}
In this paper, we present the first---to our best knowledge---example of a reduced polytope in three-dimensional Euclidean space. As the third author \cite{Polyanskii2016} has already pointed out, the existence of reduced polytopes in Euclidean space remains open starting from dimension four. Moreover, already finding a reduced polytope in three-dimensional space with different combinatorial structure than the one presented in \hyperref[chap:reduced_polytope]{Section~\ref*{chap:reduced_polytope}} seems to be a non-trivial task. Finally, it has to be checked to which amount \hyperref[thm:nonreduced]{Theorem~\ref*{thm:nonreduced}} can be generalized to higher dimensions.

\bigskip

\textbf{Acknowledgements.} We would like to thank Horst Martini for encouraging us in the search of reduced polytopes, and Ren\'{e} Brandenberg and Undine Leopold for fruitful discussions.
We also thank Alexandr Golovanov for bringing the idea to consider spherical images of polytopes to our attention. This idea helped us to construct an example of a reduced polytope.

\providecommand{\bysame}{\leavevmode\hbox to3em{\hrulefill}\thinspace}
\providecommand{\href}[2]{#2}

\end{document}